\documentclass[11pt]{amsart}

\usepackage{amsmath, amssymb, xypic}

\newcommand{\bbN}{{\mathbb N}}

\def\cK{{\mathcal K}}

\newcommand{\cZ}{{\mathcal Z}}

\newcommand{\cC}{\mathcal C}

\newcommand{\cB}{\mathcal B}

\newcommand{\bbP}{\mathbb P}
\newcommand{\e}{\varepsilon}

\newtheorem{thm}{Theorem}

\newtheorem{corollary}[thm]{Corollary}

\newtheorem{claim}[thm]{Claim}

\newtheorem{lemma}[thm]{Lemma}
\newtheorem{prop}[thm]{Proposition}

\theoremstyle{definition}

\newtheorem{definition}[thm]{Definition}

\DeclareMathOperator{\Fin}{Fin}

\newcounter{my_enumerate_counter}
\newcommand{\pushcounter}{\setcounter{my_enumerate_counter}{\value{enumi}}}
\newcommand{\popcounter}{\setcounter{enumi}{\value{my_enumerate_counter}}}

\newcommand{\FileName}[1]{\thanks{Filename: {\tt #1}}}

\newcommand{\cP}{{\mathcal P}}

\newcommand{\cU}{{\mathcal U}}
\newcommand{\cV}{\mathcal V}

\newcommand{\bp}{\mathbf p}
\newcommand{\bq}{\mathbf q}
\newcommand{\br}{\mathbf r}
\newcommand{\bs}{\mathbf s}

\newcommand{\ba}{\mathbf a}
\newcommand{\bb}{\mathbf b}
\newcommand{\bc}{\mathbf c}

\newcommand{\NneN}{\bbN^{\nearrow\bbN}}

\title[The relative commutant of separable C*-algebras]{The relative commutant of separable C*-algebras of real rank zero}

\author{Ilijas Farah}

\address{Department of Mathematics and Statistics\\
York University\\
4700 Keele Street\\
North York, Ontario\\ Canada, M3J 1P3}

\email{ifarah@mathstat.yorku.ca}

\urladdr{http://www.math.yorku.ca/$\sim$ifarah}

\thanks{Partially supported by NSERC.
I would like to thank N. Christopher Phillips for many useful
comments on the first draft of this paper. In this version
Theorem~\ref{T1} was proved only for UHF algebras, and Chris's
suggestion to use of $\leq$ instead of $\preceq$ helped me extend the
result to its present form.}

\FileName{2008i15-non-unique-commutant.tex}

\date{\today}
\begin{document}

\begin{abstract}
We answer a question of E. Kirchberg (personal
communication): does the relative commutant of a separable C*-algebra
in its ultrapower depend on the choice of the ultrafilter?
\end{abstract}

\maketitle

All algebras and all subalgebras in this note are C*-algebras and
C*-subalgebras, respectively, and all ultrafilters are nonprincipal
ultrafilters on~$\bbN$.  Our C*-terminology is standard (see e.g.,
\cite{Black:Operator}).

In the following $\cU$ ranges over nonprincipal ultrafilters on
$\bbN$. With  $A^{\cU}$ denoting the (norm, also called C*-)
ultrapower of a C*-algebra $A$ associated with  $\cU$ we have
\[
F_{\cU}(A)=A'\cap A^{\cU},
\]
the relative commutant of $A$ in its ultrapower. This invariant plays
an important role in \cite{KircPhi:Embedding} and
\cite{Kirc:Central}.

\begin{thm} \label{T1}
For every separable infinite-dimensional C*-algebra $A$ of real rank
zero the following are equivalent.
\begin{enumerate}
\item \label{T1.1} $F_{\cU}(A)\cong F_{\cV}(A)$ for any two  nonprincipal ultrafilters $\cU$ and $\cV$ on $\bbN$.
\item \label{T1.2} $A^{\cU}\cong A^{\cV}$ for any two  nonprincipal ultrafilters $\cU$ and $\cV$ on $\bbN$.
\item\label{T1.3} The Continuum Hypothesis.
\end{enumerate}
\end{thm}

The equivalence of \eqref{T1.3} and \eqref{T1.2} in Theorem~\ref{T1}
for every infinite-dimensional C*-algebra $A$ of cardinality
$2^{\aleph_0}$ that has arbitrarily long finite chains in the
Murray-von Neumann ordering of projections was proved in
\cite[Corollary~3.8]{GeHa}, using the same Dow's result from
\cite{Do:Ultrapowers} used here.

We shall prove \eqref{T1.1} implies \eqref{T1.3} and \eqref{T1.2}
implies \eqref{T1.3} in Corollary~\ref{Cor1} below. The reverse
implications are well-known consequences of countable saturatedness
of ultrapowers associated with nonprincipal ultrafilters on $\bbN$
(see \cite[Proposition~7.6]{BYBHU}). The implication from (3) to (1)
holds for every separable C*-algebra $A$ and the implication from (3)
to (2) holds for every C*-algebra $A$ of size $2^{\aleph_0}$.  The
point is that if $A$ is separable then the isomorphism between
diagonal copies of $A$ extends to an isomorphism between the
ultrapowers. Countable saturation of $A^{\cU}$ can be proved directly
from its analogue, due to Keisler, in classical model theory. This
also follows from the argument in \cite[Theorem~3.2 and
Remark~3.3]{GeHa}.

While the Continuum Hypothesis implies that any two ultrapowers of
$\cB(H)$ associated with nonprincipal ultrafilters on $\bbN$ are
isomorphic, it does not imply that the relative commutants of
$\cB(H)$ in those ultrapowers are isomorphic. As a matter of fact, it
implies the opposite (see~\cite{FaPhiSte:Relative}).

For a C*-algebra $A$ let $\cP(A)=\{p: p\in A$ is a projection$\}$
ordered by $p\leq q$ if and only if $pq=p$.
 Our proof depends on the
analysis of types of gaps in $\cP(A'\cap A^{\cU})$ (see
Definition~\ref{D1}).
 Gaps in  $\cP(\bbN)/\Fin$ and related quotient structures are
well-studied; for example, analysis of such gaps is very important in
the consistency proof of the statement `all Banach algebra
automorphisms of $C(X)$ into some Banach algebra are continuous' (see
\cite{DaWo:Introduction}). It was recently discovered that the
gap-spectrum of $\cP(\cC(H))$ (where $\cC(H)$ is the Calkin algebra,
$\cB(H)/\cK(H))$) is much richer than the gap-structure of
$\cP(\bbN)/\Fin$ (\cite{ZaaV:Gaps}).

\subsection*{Notational convention} We denote elements of
ultraproducts by boldface Roman letters such as $\bp$ and their
representing sequences by  $p(n)$, for $n\in \bbN$. We shall follow
von Neumann's convention and identify a natural number $n$ with the
set $\{0,\dots, n-1\}$. The symbol $\omega$ is used for ultrafilters
in the operator algebra literature and it is reserved for the least
infinite ordinal in the set-theoretic literature. I will avoid using
it in this note.

 By  $\sigma(a)$ we denote the spectrum of a normal operator $a$.
Lemma~\ref{L.-1} below is well-known. A sharper result can be found
e.g.,  in \cite[Lemma~2.5.4]{Lin:Introduction} but  we include
 a proof for reader's convenience.

\begin{lemma} \label{L.-1} For a self-adjoint  $a$ and a projection $r$, if $\|a-r\|<\e<1$ then
$\sigma(a)\subseteq (-2\sqrt\e,2\sqrt\e)\cup
(1-2\sqrt\e,1+2\sqrt\e)$. If in addition $\e<1/16$ then there is a
projection $r'$ in $C^*(a)$ such that $\|r'-a\|<2\sqrt \e$.
\end{lemma}

\begin{proof} Since $\|a\|<1+\e<2$, we have $\|a^2-a\|\leq \|a(a-r)\|+\|r(a-r)\|+\|a-r\|<4\e$.
Thus $|x(1-x)|<4\e$ for all $x\in \sigma(a)$ and in turn
$|x|<2\sqrt\e$ or $|1-x|<2\sqrt \e$.

Now assume $\e<1/16$. In this case $1/2\notin \sigma(a)$. Define a
continuous function $f$ with domain $\sigma(a)$ as follows. Let
$f(t)=0$ for $-\infty<t<1/2$ and $f(t)=1$ for $1/2\leq t<\infty$.
Since $|f(t)-t|<2\sqrt \e$ for all $t\in \sigma(a)$,  $f(a)$ is a
projection in $C^*(a)$ as required. \end{proof}

A representing sequence $p(n)$ of a projection $\bp$ in an ultrapower
can be chosen so that each $p(n)$ is a projection (see
\cite[Proposition~2.5 (1)]{GeHa}, this also follows immediately from
\cite[Lemma~4.2.2]{Lor:Lifting} or
\cite[Lemma~2.5.5]{Lin:Introduction}).

\begin{lemma} \label{L4}
For projections $\bp,\bq$ in $A^{\cU}$ the following are equivalent.
\begin{enumerate}
\item \label{L4.1} $\bp\leq \bq$,
\item \label{L4.2} There is a representing sequence $p'(i)$, for $i\in \bbN$, of $\bp$ such
that $p'(i)\leq q(i)$ for all $i$.
\item \label{L4.3} There is a representing sequence $q'(i)$, for $i\in \bbN$, of $\bq$ such
that $p(i)\leq q'(i)$ for all $i$.
\end{enumerate}
\end{lemma}

\begin{proof} Both  \eqref{L4.3} implies \eqref{L4.1} and
\eqref{L4.2} implies \eqref{L4.1} are trivial. We shall prove
\eqref{L4.1} implies \eqref{L4.2}.   Assume $\bp\leq \bq$. For every
$n\geq 1$ the set
\[
X_n=\{j: \|q(j)p(j)q(j)-p(j)\|<1/(4n)\}
\]
belongs to  $\cU$. We may assume $\bigcap_n X_n=\emptyset$. Let
$p'(j)=0$ if $j\notin X_0$. If $j\in X_n\setminus X_{n+1}$ then
Lemma~\ref{L.-1}, with $a(j)=q(j)p(j)q(j)$, implies there is a
projection $p'(j)\in C^*(a(j))$ such that $\|p'(j)-a(j)\|<1/(2\sqrt
n)$. Then $p'(j)\leq q(j)$ and $\|p'(j)-p(j)\|<1/\sqrt n$ for all
$j\in X_n$. Therefore $p'(j)$, for $j\in \bbN$, is a representing
sequence of $\bp$ as required.

In order to prove \eqref{L4.1} implies \eqref{L4.3} apply the above
to $1-\bp\geq 1-\bq$ in the ultrapower of the unitization of $A$ to
find an appropriate representing sequence for $1-\bq$.
\end{proof}

By $\NneN$ we denote the set of all nondecreasing functions $f$ from
$\bbN$ to $\bbN$ such that $\lim_n f(n)=\infty$, ordered pointwise.
Write $f\leq_{\cU} g$ if $\{n: f(n)\leq g(n)\}\in \cU$ and denote the
quotient linear ordering by  $\NneN/\cU$.

Following \cite{Do:Ultrapowers}, for an ultrafilter $\cU$ we write
$\kappa(\cU)$ for the \emph{coinitiality} of $\NneN/\cU$, i.e., the
minimal cardinality of  $X\subseteq \NneN$ such that for every $g\in
\NneN$ there is $f\in X$ such that $f\leq_{\cU} g$. (It is not
difficult to see that this is equal to $\kappa(\cU)$ as defined in
 \cite[Definition~1.3]{Do:Ultrapowers}.)

\begin{definition}\label{D1}
Let $\lambda$ be a cardinal. An \emph{$(\aleph_0,\lambda)$-gap} in a
partially ordered set $\bbP$ is a pair consisting of a
$\leq_{\bbP}$-increasing family $\ba_m$, for $m\in \bbN$, and a
$\leq_{\bbP}$-decreasing family $\bb_\gamma$, for $\gamma<\lambda$,
such that $\ba_m\leq_{\bbP} \bb_\gamma$ for all $m$ and $\gamma$ but
there is no $\bc\in \bbP$ such that $\ba_m\leq_{\bbP} \bc$ for all
$m$ and $\bc\leq_{\bbP}\bb_\gamma$ for all $\gamma$.
\end{definition}

Assume $r^0(n)\leq r^1(n)\leq \dots \leq r^{l(n)-1}(n)$ are
projections in $A$ and $\lim_{n\to \infty}l(n)=\infty$. For $h\colon
\bbN\to \bbN$ define $\br^h$ via its representing sequence (let
$r^i(n)=r^{l(n)-1}(n)$ for $i\geq l(n)$)
\[
r^h(n)=r^{h(n)}(n).
\]
Let $\bp_m=\br^{\bar m}$, where $\bar m(j)=m$ for all $j$.

\begin{lemma}\label{L5}
With notation from the previous paragraph, for every projection $\bs$
in $A^{\cU}$ such that $\bp_m\leq \bs$ for all
 $m$ there is $h\colon \bbN\to \bbN$ such that $\bp_m\leq \br^h$ for all $m$
 and $\br^h\leq
 \bs$.
\end{lemma}

\begin{proof} Since $\bp_m\leq \bs$, for each $m\in \bbN$ the set
\[
X_m=\left\{i: \|r^m(i)s(i)-r^m(i)\|<1/m\right\}
\]
belongs to $\cU$. Since the value of $\|r^m(i)s(i)-r^m(i)\|$ is
increasing in $m$ we have $X_m\supseteq X_{m+1}$. We may assume
$\bigcap_m X_m=\emptyset$. Define $h\colon \bbN\to \bbN$ by letting
$h(i)=0$ for $i\notin X_0$ and for $i\in X_m\setminus X_{m+1}$ let
$h(i)=m$.

For each $m$ and $i\in X_m$ we have $h(i)\geq m$ and therefore
$\br^h\geq \bp_m$. Also, $i\in X_m$ implies
$\|r^h(i)s(i)-r^h(i)\|<1/m$ hence $\br^h\leq \bs$.
\end{proof}

 The proof of Proposition~\ref{P0} was inspired by Alan
Dow's \cite[Proposition~1.4]{Do:Ultrapowers}. Dow's result was
independently proved by Saharon Shelah and can be found in
\cite{She:Classification}.\footnote{I could not find it, but it
should be somewhere in Chapter VI.}

By $A_{\leq 1}$ we denote the unit ball of a C*-algebra $A$.

\begin{prop}\label{P0} Assume $A$ is a separable C*-algebra and there are finite
self-adjoint sets $F_0\subseteq F_1\subseteq F_2\subseteq \dots
\subseteq A_{\leq 1}$ whose union is dense in $A_{\leq 1}$ and such
that for each $n$ there is a $\leq$-increasing chain $\cC_n$ of
projections in $B_n=F_n'\cap A$ of length at least $n$.

Then for every nonprincipal ultrafilter $\cU$ on $\bbN$ and every
cardinal $\lambda$ there is an $(\aleph_0,\lambda)$-gap in
$\cP(A'\cap A^{\cU})$ if and only if $\kappa(\cU)=\lambda$.
\end{prop}

\begin{proof} First we prove the converse implication.
Assume $g_\gamma$, for $\gamma<\lambda=\kappa(\cU)$, is a
$\leq_{\cU}$-decreasing and $\leq_{\cU}$-unbounded below chain of
functions in $\NneN$. Let $0=r^0(n)\leq r^1(n)\leq \dots \leq
r^{n-1}(n)$ be an enumeration of $\cC_n$.

\begin{claim} \label{C0} For all $f,g$ in $\NneN$ the following are
equivalent.
\begin{enumerate}
\item  \label{C0.I1} $f\leq _{\cU} g$,
\item  \label{C0.I2}  $\br^f\leq\br^g$,
\end{enumerate}
\end{claim}

\begin{proof} Assume $f\leq_{\cU} g$. Then $X=\{j: f(j)\leq g(j)\}\in
\cU$ and $\br^f(j)\leq \br^g(j)$ for all $j\in X$ hence \eqref{C0.I2}
follows.   If  $f\not\leq_{\cU} g$ then $X=\{j: f(j)> g(j)\}\in \cU$
and for all $j\in X$ we have $\|r^f(i)r^g(i)-r^g(i)\|=1$, hence
$\br^f\nleq\br^g$.
\end{proof}

Let $\bq_\gamma=\br^{g_\gamma}$, for $\gamma<\lambda$.  By
Claim~\ref{C0} we have
\[
\bp_m\leq \bp_{m+1}\leq \bq_\delta\leq \bq_\gamma
\]
for all $m$ and all $\gamma<\delta<\lambda$. All of $\bp_m$ and
$\bq_\gamma$ belong to $A'\cap A^{\cU}$.

We shall show that this family forms a gap in $\cP(A^{\cU})$ (and
therefore it forms a gap in $\cP(A'\cap A^{\cU})$. Assume $\bs\in
A^{\cU}$ is such that $\bs\leq \bq_\gamma$ for all $\gamma$. By
Lemma~\ref{L5} there is $h$ such that $\bp_m\leq \br^h\leq \bs$ for
all $m$. By Claim~\ref{C0} we have $h\leq_{\cU} g_\gamma$ for all
$\gamma$ and $\bar m\leq_{\cU} h$ for all $m$, a contradiction.

In order to prove the direct implication, assume that $\bp_m$,
$\bq_\gamma$ form an $(\aleph_0,\lambda)$-gap in $\cP(A'\cap
A^{\cU})$. By successively using Lemma~\ref{L4} for $m=1,2,\dots$
find representing sequences $p_m(i)_{i\in \bbN}$,  for $\bp_m$ such
that $p_m(i)\leq p_{m+1}(i)$ for all $i$. Choose an increasing
sequence $0=m_0<m_1<m_2<\dots$ such that the following holds for all
$k$.
\begin{enumerate}
\item  [(*)]for all $j<m_k$ and all $a\in F_{m_k}$, if $l\geq m_{k+1}$
then $\|[p_j(l),a]\|<1/k$.
\end{enumerate}
For $n\in \bbN$ and $i$ such that for some $k$ we have $i<m_k$ and
$m_{k+1}\leq n$ let $r^i(n)=p_i(n)$. Thus we have projections
\[
r^0(n)\leq r^1(n)\leq\dots\leq r^{m_k}(n)
\]
whenever $n\geq m_{k+1}$. For $h\colon \bbN\to \bbN$ define $\br^h$
as in the paragraph before Lemma~\ref{L5}, by its representing
sequence (let $r^i(n)=r^{m_k}(n)$ if $i\geq m_k$)
\[
r^h(n)=r^{h(n)}(n).
\]

\begin{claim}\label{C0.1}
If $h\colon \bbN\to \bbN$ then $\br^h\in A'\cap A^{\cU}$.
\end{claim}

\begin{proof} Fix any $b$ in the unit ball of $A$ and $\e>0$. If $k>1/\e$
and there is $b'\in F_{2k}$ satisfying $\|b-b'\|<\e/2$
then for $i>n_{2k}$ in $Y$ we have that $\|[p_j(i),b']\|<\e/2$ and
therefore $\|[\br^h(i),b]\|<\e$ for $\cU$-many $i$.
\end{proof}

 Using Lemma~\ref{L5} for each $\bq_\gamma$ find $h_\gamma$
such that  $\br^\gamma=\br^{h_\gamma}$ satisfies $\bp_i\leq
\br^\gamma\leq \bq_\gamma$ for all $i$. Since $\NneN/\cU$ is a linear
ordering and $\lambda$ is a regular cardinal, we can find a cofinal
subset $\cZ$ of $\lambda$ such that for $\gamma<\delta$ in $\cZ$ we
have $\br^\delta\leq \br^\gamma$. By reenumerating we may assume
$\cZ=\lambda$ and then $\br^\gamma$, for $\gamma\in \cZ$, together
with $\bp_i$, for $i\in \bbN$, form an $(\aleph_0,\lambda)$-gap.
However, $\br^\delta\leq\br^\gamma$ is equivalent to
$h_\delta\leq_{\cU} h_\gamma$, and therefore $h_\gamma$, for
$\gamma<\lambda$, form a $\leq_{\cU}$-decreasing and
$\leq_{\cU}$-unbounded below sequence in $\NneN/\cU$, and therefore
$\lambda=\kappa(\cU)$.
\end{proof}

The proof of Proposition~\ref{P0} can be modified (by
removing some of its parts) to a proof of the following.

\begin{prop}\label{P0-} Assume $A$ is a separable C*-algebra and
 $\cP(A)$ has arbitrarily long finite chains. Then for every
nonprincipal ultrafilter $\cU$ on $\bbN$ and every cardinal $\lambda$
there is an $(\aleph_0,\lambda)$-gap in $\cP(A^{\cU})$ if and only if
$\kappa(\cU)=\lambda$. \qed
\end{prop}

\begin{corollary} \label{Cor1} Assume the Continuum Hypothesis fails.
If $A$ is an infinite-dimensional separable  C*-algebra of real rank
zero then there are nonprincipal ultrafilters $\cU$ and $\cV$ on
$\bbN$ such that $F_{\cU}(A)\not\cong F_{\cV}(A)$ and
$A^{\cU}\not\cong A^{\cV}$.
\end{corollary}

\begin{proof}  By
\cite[Theorem~2.2]{Do:Ultrapowers} we can find $\cU$ and $\cV$ so
that $\kappa(\cU)=\aleph_1$ and $\kappa(\cV)=\aleph_2$ (here
$\aleph_1$ and $\aleph_2$ are the least two uncountable cardinals;
all that matters for us is that they are both less or equal than
$2^{\aleph_0}$ and different). Therefore $\cP(A'\cap A^{\cU})$ has an
$(\aleph_0,\aleph_1)$-gap while $\cP(A'\cap A^{\cV})$ does not, and
$A'\cap A^{\cU}$ and $A'\cap A^{\cV}$ cannot be isomorphic.

It remains to prove that if $A$ is an infinite-dimensional C*-algebra
of real rank zero then $\cP(A)$ has an infinite chain of projections.
We may assume $A$ is unital.  Recursively find a decreasing sequence
$r_n$ for $n\in \bbN$ in $\cP(A)$ so that $r_nAr_n$ is
infinite-dimensional for all $n$. Assume $r_n$ has been chosen. Since
$A$ has real rank zero, in $r_n A r_n$ we can fix a projection
$q\notin \{0,r_n\}$. If $q A_n q$ is infinite-dimensional then let
$r_{n+1}=q$. Otherwise, let $r_{n+1}=r_n-q$ and note that $r_{n+1} A
r_{n+1}$ is infinite-dimensional.
\end{proof}

It is likely that Theorem~\ref{T1} and Corollary~\ref{Cor1} can be
extended to all infinite-dimensional separable C*-algebras (possibly
by considering the Cuntz ordering of positive elements instead
of~$\cP(A)$).

\providecommand{\bysame}{\leavevmode\hbox
to3em{\hrulefill}\thinspace}
\providecommand{\MR}{\relax\ifhmode\unskip\space\fi MR }
\providecommand{\MRhref}[2]{%
  \href{http://www.ams.org/mathscinet-getitem?mr=#1}{#2}
} \providecommand{\href}[2]{#2}

\end{document}